\documentclass[11pt]{imsart}
\RequirePackage[OT1]{fontenc}
\RequirePackage{amsthm,amsmath}
\RequirePackage[numbers]{natbib}
\RequirePackage[colorlinks,citecolor=blue,urlcolor=blue]{hyperref}

\usepackage{graphicx,amsmath,amsthm,amsfonts,amssymb}
\usepackage{booktabs, rotating, float}
\usepackage[hmargin=3cm,vmargin=3cm]{geometry}
\usepackage{mathtools}
\usepackage{bbm}
\mathtoolsset{showonlyrefs}
\usepackage{enumerate}
\usepackage{color}
\usepackage{listings}
\usepackage{epstopdf}
\usepackage{makecell}
\usepackage{tikz}
\usetikzlibrary{decorations.pathreplacing}
\usepackage{array}
\usepackage{caption}
\usepackage{subcaption}
\usepackage{soul}
\usepackage{hyperref}
\usepackage{verbatim}
\usepackage{siunitx}
\usetikzlibrary{calc} \tikzset{>=latex}
\usetikzlibrary{calc,decorations.markings}
\usetikzlibrary{positioning,shapes,decorations.text}
\theoremstyle{plain}
\newtheorem{thm}{Theorem}[section]
\newtheorem{cor}[thm]{Corollary}

\newtheorem{prop}[thm]{Proposition}
\newtheorem{lemma}[thm]{Lemma}
\newtheorem{assumption}[thm]{Assumption}
\theoremstyle{definition}
\newtheorem{defn}[thm]{Definition}

\newtheorem{eg}[thm]{Example}
\newtheorem{fact}[thm]{Fact}
\newtheorem{observe}[thm]{Observation}

\numberwithin{equation}{section}

\newcommand{\rpm}{\sbox0{$1$}\sbox2{$\scriptstyle\pm$}
  \raise\dimexpr(\ht0-\ht2)/2\relax\box2 }

\tikzstyle{nd} = [anchor=base, inner sep=0pt]
\tikzstyle{ndpic} = [remember picture, baseline, every node/.style={nd}]

\def\beq{\begin{equation}}
\def\eeq{\end{equation}}
\def\ba{\begin{enumerate}[(a)]}
\def\bei{\begin{enumerate}[(i)]}
\def\be{\begin{enumerate}[(1)]}
\def\ee{\end{enumerate}}
\def\bi{\begin{itemize}}
\def\ei{\end{itemize}}
\def\beg{\begin{eg}}
\def\eeg{\end{eg}}
\def\bd{\begin{defn}}
\def\ed{\end{defn}}
\def\bt{\begin{thm}}
\def\et{\end{thm}}
\def\bl{\begin{lemma}}
\def\el{\end{lemma}}
\def\bfac{\begin{fact}}
\def\efac{\end{fact}}

\def\bc{\begin{cor}}
\def\ec{\end{cor}}
\def\bp{\begin{prop}}
\def\ep{\end{prop}}
\def\bo{\begin{observe}}
\def\eo{\end{observe}}
\def\bas{\begin{assumption}}
\def\eas{\end{assumption}}
\def\RR{\mathbb{R}}

\def\ZZ{\mathbb{Z}}
\def\NN{\mathbb{N}}

\def\ii{\item}

\def\beg{\begin{eg}}
\def\eeg{\end{eg}}

\makeindex

\numberwithin{equation}{section}
\numberwithin{table}{section}

\startlocaldefs
\endlocaldefs

\begin{document}

\begin{frontmatter}

\title{Moments of the SHE under delta initial measure}
\runtitle{Moments of SHE.}
\runauthor{Ghosal}

\begin{aug}
  \author{\fnms{Promit} \snm{Ghosal}\ead[label=e1]{pg2475@columbia.edu}}
%
%
%
%
%
%
\address{ Department of Statistics, Columbia University, 1255 Amsterdam Avenue, New York, NY 10027\\ \printead{e1}}
\end{aug}

\begin{abstract}
We give a rigorous proof of the contour integral formulas of the one-point moments of the stochastic heat equation (SHE) started from the delta initial measure at the origin. These formulas were conjectured in \cite{BC14} (see also \cite{DR02, DOT03}).  Our proof is based on a correspondence between the SHE and the Airy point process which was proved in \cite[Theorem~1]{BorGor16} using the formula of \cite[Theorem~1.1]{Amir11}. 
\end{abstract}

%
%

\end{frontmatter}









\section{Introduction}\label{Intro}
The stochastic heat equation with delta initial measure at the origin is given as 
\begin{align}\label{eq:SHEDef}
\frac{\partial }{\partial T}\mathcal{Z}(T,X) &= \frac{1}{2} \frac{\partial^2}{\partial^2 X} \mathcal{Z}(T,X) +\frac{1}{2} \mathcal{Z}(T,X) \mathcal{W}(T,X) \qquad X\in \RR, t\in \RR_{+}\\
\mathcal{Z}(0,X) & = \delta_{X=0}. 
\end{align} 
Here, $\mathcal{W}(T,X)$ is space time white noise which is a distribution valued Gaussian random field with the following correlation structure 
\begin{align}\label{eq:WNCorr}
\mathbb{E}\Big[\mathcal{W}(T_1,X_1)\mathcal{W}(T_2, X_2)\Big] = \delta_{T_1=T_2} \delta_{X_1= X_2}.
\end{align}
  For the solution theory of the SHE, we refer to \cite{Walsh86,Corwin12,Quastel12}. Applications of the SHE range from modeling the density of the particles diffusing under the space-time random environment or random drifts \cite{Mol96, BarraquandCorwin15, CG17} to the contnuous directed random polymers \cite{AKQ1}. The logarithm of the SHE is called the \emph{Cole-Hopf} solution of the \emph{Kardar-Parisi-Zhang} (KPZ) equation which is a prototype for random growth processes and a testing ground for the study of nonlinear stochastic PDEs. Moments of the SHE are important ingredients in extracting information about intermittency \cite{Amir11, Bertini1995, BC14, CK10, CK12, CJKS}, regularity \cite{SS00, SS02}, tail decay \cite{CDR10, CJK} etc. In this article, we prove one-point moment formulas of $\mathcal{Z}(T,X)$. Our main result is stated as follows.
\bt\label{MainTheo} 
Let $\mathcal{Z}(T,X)$ be the unique solution of \eqref{eq:SHEDef}. Then, for any $k\in \NN$
\begin{align}\label{eq:SHEMom}
\mathbb{E}_{\mathrm{SHE}}\left[(\mathcal{Z}(T, X))^{k}\right] = \frac{1}{(2\pi \mathbf{i})^{k}} \int_{\mathcal{C}_1}\ldots \int_{\mathcal{C}_k} \prod_{1\leq A < B\leq k} \frac{z_{A} - z_{B}}{z_{A} - z_{B}-1} e^{\frac{T}{2}\sum_{j=1}^{k} z^2_j + X\sum_{j=1}^{k}z_j} \prod_{j=1}^{k} dz_j  &&&&
\end{align}
where $\mathcal{C}_j$ is the line $\alpha_j +\mathbf{i}\RR$ such that $\alpha_1>\alpha_2+2> \ldots >\alpha_k+(k-1)$
\et   

Recently, \cite[Theorem~2]{BorGor16} established an identity between the integer moments of $\mathcal{Z}(T,0)$ and the complete homogeneous symmetric functional moments of the \emph{Airy point process}. Their proof relies on the moment formulas of \cite[Proposition~5.4.8]{BC14} which are same as in Theorem~\ref{MainTheo}. Here, we demonstrate an alternative way deriving those formulas. The main technical innovation of this paper is to prove the correspondence (see Proposition~\ref{AiryToSHE}) between the moments of $\mathcal{Z}(T,0)$ and the Airy point process without using the formula of \cite[Proposition~5.4.8]{BC14}. In combination with a rigorous derivation of the Airy moments (see Lemma~\ref{ExpKSumLem}), this correspondence will provide a rigorous proof of Theorem~\ref{MainTheo}. Here, we would like to stress that the derivation of the Airy moments were known from the work of \cite{BorGor16}. For the sake of completeness, we will give a brief outline of their proof avoiding much details. In \cite[Proposition~6.2.3]{BC14}, the authors showed that the right hand side of \eqref{eq:SHEMom} solves the delta-Bose gas with delta potential. Combining our main theorem with \cite[Proposition~6.2.3]{BC14} yields the proof of the conjecture that the one-point moments of $\mathcal{Z}(T,X)$ coincides with the solution of the delta-Bose gas with delta potential. However, the uniqueness of the solution is not clear yet.  

 As a main tool, we use identity $(1)$ of \cite[Theorem~1]{BorGor16} which relates the Laplace transform of $\mathcal{Z}(T,0)$ with some multiplicative functional of the Airy point process. More specifically, they showed 
  \begin{align}\label{eq:identity}
  \mathbb{E}_{\mathrm{Airy}}\Big[\prod_{p=1}^{\infty}\frac{1}{1+u\exp(C\mathbf{a}_p)}\Big]=\mathbb{E}_{\mathrm{SHE}}\Big[\exp\big(-u\mathcal{Z}(T,0)\exp(T/24)\big)\Big] , \quad \forall u\in \RR_{+}
  \end{align}
  where $C=(T/2)^{1/3}$ and $\mathbf{a}_1\geq \mathbf{a}_2\geq  \ldots $ are the ordered points of the Airy point process (see Section~\ref{SHEMoment} for its definition). This identity comes from some determinantal manipulation of the formulas in \cite[Theorem~1.1]{Amir11} and does not require the moment formula of $\mathcal{Z}(T,X)$ as an input. One should note that similar identities have also been observed in other models, namely, between the stochastic higher spin vertex model and the Macdonald measure \cite[Theorem~4.2]{B16Aug}, the asymmetric simple exclusion process (ASEP) and the discrete Laguerre ensemble \cite[Theorem]{BO16} etc. In fact, \eqref{eq:identity} follows from those identities by taking appropriate limits. From Taylor expanding the right hand side of \eqref{eq:identity}, one sees that coefficient of $(-u)^{k}$ is the $k$th moment of $\mathcal{Z}(T,0)$ upto some constant. However, the series obtained after Taylor expanding the the left hand side of \eqref{eq:identity} is divergent in nature. Because of this, one faces difficulties in interchanging derivatives and expectation on the right hand side \eqref{eq:identity}. We overcome these difficulties by considering only finitely many terms in the expansion of $[1+u\exp(C\mathbf{a}_p)]^{-1}$ for $p=1,2, \ldots $ and bounding the contributions of the remainder terms (see Section~\ref{ASHE} for more details) appropriately.

     The moment formulas of \eqref{eq:SHEMom} were formally derived in \cite{BC14} in two different ways. One of those two ways hinges upon the folklore that the moments of $\mathcal{Z}(T,X)$ solves the attractive \emph{delta-Bose gas}. In \cite[Proposition~6.2.3]{BC14}, the right hand side of \eqref{eq:SHEMom} is shown to be a solution of the  attractive delta-Bose gas with delta potential. Under the assumption that the folklore is true, their result indicates the moment formulas of Theorem~\ref{MainTheo}. However, a rigorous proof of this claim is missing except when $k=2$ (see \cite{AGHH}). Proof for the case $k=2$ relies on the observation that the moments of the SHE with smoothed out (in space) white noise solve the Lieb-Liniger many body problem with smoothed delta potential. It then suffices to show that the smoothed moments converges to the moments without smoothing and likewise, the solution of the Lieb-Liniger system with smooth potential converges to the solution of the delta-Bose gas with delta interaction. However, showing this twofold convergence is nontrivial and becomes complicated for higher order moments. To our knowledge, there is no follow up work where this approach is made to work for the moments of order greater than $2$. 
  
  The second approach of \cite{BC14} is based on the observation that the discrete and semi-discrete directed polymers converges in distribution to the SHE under the intermediate disorder scaling. See \cite{alberts2014, DR02} for the discrete case and \cite{Nica16} for the semi-discrete case. Moments of the semi-discrete polymer measure are given by some integral formulas of the Whittaker measure which is a degeneration of the Macdonald measure. By taking limit of the integral formulas of the Whittaker measure, \cite[Section~5.4.2]{BC14} obtained the multipoint moment formulas of the SHE. However, the convergence in distribution is not good enough to support the convergence of moments. To show the moments also converge, one needs uniform tail bounds of the pre-limiting object.
   Prof. Vadim Gorin kindly informed us that it is indeed possible to obtain such bounds via Markov's inequality and analysis of known contour integral formulas for the semi-discrete polymer. Combining these tail bounds with the known weak convergence provides an alternative derivation of \eqref{eq:SHEMom}. For the sake of completeness, we include some details of this approach in Section~\ref{AlrMethod}. We should also point out that using this approach, it is possible to to obtain multi-point moment formulas of the SHE as conjectured in \cite[Proposition~5.4.6]{BC14}. For more details, we refer to \cite[Corollary~1.14]{Nica16}.

There is yet another approach to derive the moment formulas of Theorem~\ref{MainTheo}\footnote{Communicated to us by Guillaume Barraquand}. The main idea of this approach relies upon an identity between the distribution of $\mathcal{Z}(0,T)$ and a random infinite series involving the Airy point process and an infinite sequence of independent chi-squared  random variables. This identity is proved in \cite[Theorem~1.3]{GS18} and is stated as follows 
 \begin{align}\label{eq:RandomSeries}
 \mathcal{Z}(T,0)\exp(T/24) \stackrel{d}{=} 2\sum_{p=1}^{\infty} u^2_p\exp(C\mathbf{a}_p)
\end{align}    
  where $\{u^2_p\}_{p\in \NN}$ is a sequence of independent 
 chi-squared random variables with $1$ degree of freedom. Hence, the $k$-th moment of $\mathcal{Z}(0,T)$ is given by the $k$-th moment of the random infinite series of the right hand side of \eqref{eq:RandomSeries} modulo the constant $\exp(-kT/24)$. Now, the contour integral formula of Theorem~\ref{MainTheo} can be obtained by expanding the $k$-th moments of the random infinite series and applying the moment formulas of the chi-squared random variables. However, we do not pursue the details of these computations in this paper.

  Recently, identity \eqref{eq:identity} played a  key role in \cite{CG18} to find the crossover behavior of the lower tail probability of the Cole-Hopf solution of the KPZ equation. In \cite{CGKLT18}, this identity was at the focal point in deriving the lower tail large deviation rate function of the KPZ equation. In an upcoming work \cite{CG18b}, the main result of this paper will be used to estimate the upper tail fluctuation of the KPZ equation. We should also point out that \cite[Proposition~2.3]{Bertini1995} obtained some moment formulas of the SHE under the general initial conditions. It is not clear whether those moment formulas (even for $k=2$) coincide with the right hand side of \eqref{eq:SHEMom} under the delta initial measure. Those formulas of \cite{Bertini1995} were later used (for instance in \cite{CJK, ChD15}) to get the upper tail estimates of the SHE. However, for our purpose in \cite{CG18b}, we need more refined estimates which we obtain from \eqref{eq:SHEMom}. This was our primary motivation to give a rigorous proof of Theorem~\ref{MainTheo}.

The rest of the paper is organized as follows. We prove Theorem~\ref{MainTheo} in Section~\ref{ProofOfMT}. In Section~\ref{AlrMethod}, we demonstrate one of the other alternative ways of proving Theorem~\ref{MainTheo}.

\section{Proof of Theorem~\ref{MainTheo}}\label{ProofOfMT}

 Substituting $\tilde{z}_j = z_j+ \frac{X}{T}$ on the right hand side of \eqref{eq:SHEMom}, we get
\begin{align}
\text{r.h.s of \eqref{eq:SHEMom}} = \frac{e^{-\frac{kX^2}{2T}}}{(2\pi \mathbf{i})^{k}} \int_{\mathcal{C}_1}\ldots \int_{\mathcal{C}_k} \prod_{1\leq A<B \leq k} \frac{\tilde{z}_A - \tilde{z}_{B}}{\tilde{z}_A - \tilde{z}_{B}-1} e^{\frac{T}{2}\sum_{j=1}^{k}\tilde{z}^2_j} \prod_{j=1}^{k} d\tilde{z}_j.\label{eq:Subs}
\end{align}
 Using spatial stationarity of the process $\mathcal{Z}(T,X)\exp\big(\frac{X^2}{2T}\big)$ (see \cite[Theorem~1.4]{Amir11}, \cite[Proposition~1.17(1)]{CorHam16}), it suffices to show  
 \begin{align}\label{eq:ObsMomeq}
 \mathbb{E}_{\mathrm{SHE}}\left[(\mathcal{Z}(T, 0))^{k}\right] = \frac{1}{(2\pi \mathbf{i})^{k}} \int_{\mathcal{C}_1}\ldots \int_{\mathcal{C}_k} \prod_{1\leq A < B\leq k} \frac{z_{A} - z_{B}}{z_{A} - z_{B}-1} e^{\frac{T}{2}\sum_{j=1}^{k} z^2_j} \prod_{j=1}^{k} dz_j.  
 \end{align}
By residue expansion, \cite[Proposition~3.2.1]{BC14} (see also \cite[Proposition~5.1]{BBC16}) gave an alternative form of the right hand side of \eqref{eq:ObsMomeq}. Our next result shows that the moments of $\mathcal{Z}(T,0)$ matches with their alternative formula.

\bt\label{SHEMoment} 
Continuing with the notation of Theorem~\ref{MainTheo}, one has   
\begin{align}
\mathbb{E}_{\mathrm{SHE}}\Big[(\mathcal{Z}(T,0))^k\Big] = \sum_{\substack{\lambda\vdash k\\ \lambda = 1^{m_1} 2^{m_2}\ldots }}  &\frac{k!}{m_1!m_2!\ldots } \int^{\mathbf{i}\infty}_{-\mathbf{i}\infty} \frac{dw_1}{2\pi \mathbf{i}}\ldots \int^{\mathbf{i}\infty}_{-\mathbf{i}\infty} \frac{dw_{\ell(\lambda)}}{2\pi \mathbf{i}} \mathrm{det}\left[\frac{1}{w_i+\lambda_i- w_j}\right]^{\ell(\lambda)}_{i,j=1}\\
&\times \prod_{j=1}^{\ell(\lambda)}\exp\left(\frac{T}{2}\left[w^2_j + (w_j+1)^2+ \ldots (w_j+\lambda_j-1)^2\right]\right). \label{eq:Moment}
\end{align}
Here, $\lambda \vdash k$ denotes a partition $k$. For any partition  $\lambda =(\lambda_1 >\lambda_2>\ldots)$, the number of nonzero elements of $\lambda$, i.e., $|\{k: \lambda_k\neq 0\}|$ is called the length of the partition $\lambda$  and denoted by $\ell(\lambda)$. 
\et


 \textsc{Final step of proof of Theorem~\ref{MainTheo}:}
  Combining \eqref{eq:Moment} with \cite[Proposition~3.2.1]{BC14} yields \eqref{eq:SHEMom}. This completes the proof.

\subsection{Proof of Theorem~\ref{SHEMoment}}

We prove Theorem~\ref{SHEMoment} using identity \eqref{eq:identity}. In order to do so, our first step is to give an alternative derivation of the equivalence between the moments of $\mathcal{Z}(T,X)$ and the complete homogeneous symmetric functional moments of the Airy point process (see Proposition~\ref{AiryToSHE}). This equivalence was first proved in \cite[Theorem~1.2]{BorGor16} subjected to the fact that the moments $\mathcal{Z}(T,X)$ solves the delta-Bose gas with delta potential. Unlike \cite{BorGor16}, we take the route of Taylor expanding both sides of  \eqref{eq:identity} with respect $u$. At this point, it is not hard to guess that equating the coefficients of $(-u)^{k}$ on both sides of \eqref{eq:identity} one may get the desired identity of \cite[Theorem~1.2]{BorGor16}. Although, it sounds very straightforward, the exact details is convoluted due to the divergent nature of the power series obtained after Taylor's expansion. The main difficulty lies in exchanging the  derivatives and the expectation which we deal in Section~\ref{AiryToSHE}. Our second step is to show that the right hand side of \eqref{eq:Moment} coincides with the complete homogeneous symmetric polynomial moment (of order $k$) of the Airy point process upto the factor $\exp\big(\frac{kT}{24}\big)$. This is proved in Lemma~\ref{ExpKSumLem}. We state our results after a brief overview of the Airy point process.

\bd
 \emph{Airy point process} is a \emph{determinantal point process} on $\RR$. A \emph{point process} on $\RR$ is an integer-valued measure on the point configurations of the real line. Any point process $\chi$ is characterized by its correlation functions $\{\rho^{\chi}_k\}_{k\in \NN}$. For any $k\in \NN$, $\rho^{\chi}_k:\RR^k\to \RR_{\geq 0}$ is a locally integrable function such that for any  Borel sets $B_1, \ldots ,B_k\in \mathcal{B}(\RR)$ 
 \begin{align}\label{eq:Expectation}
 \mathbb{E}\Big[\prod_{i=1}^{k} \mathbbm{1}_{\chi}(B_i)\Big]= \int_{\prod_{i=1}^kB_i} \rho_k(x_1, \ldots ,x_k) dx_1 \cdots dx_k
\end{align}   
A point process $\chi$ is called \emph{determinantal} when $\rho^{\chi}(x_1, \ldots ,x_k)= \mathrm{det}[K_{\chi}(x_i,x_j)]_{1\leq i\leq j\leq k}$ for some $K_{\chi}:\RR^2\to \RR$, namely, the kernel of the point process $\chi$. The kernel of the Airy point process is given as 
\begin{align}\label{eq:Airy}
K_{\mathrm{Ai}}(x,y) = \int^{\infty}_{-\infty}\mathrm{Ai}(x+t)\mathrm{Ai}(y+t) dt.
\end{align}
\ed

\bp\label{AiryToSHE} 
Fix $T>0$ and set $C= (T/2)^{\frac{1}{3}}$. Let $\mathbf{a}_1 \geq \mathbf{a}_2 \geq \ldots $ be the ordered points of the Airy point process. Let $h_k(x_1, x_2, \ldots )= \sum_{i_1\leq i_2\leq \ldots \leq i_k}x_{i_1}x_{i_2}\ldots x_{i_k}$ be the complete homogeneous symmetric functions in variables $x_1, x_2, \ldots $. Then, for any $k\in \NN$, 
\begin{align}\label{eq:kMomEq}
 \mathbb{E}_{\mathrm{Airy}}\big[h_k(\exp(C\mathbf{a}_1),\exp(C\mathbf{a}_2), \ldots )\big] = \mathbb{E}_{\mathrm{SHE}}\Big[ \frac{(\mathcal{Z}(T,0))^k}{k!}\exp\big(\frac{kT}{24}\big)\Big]
 \end{align}
\ep

Proof of Proposition~\ref{AiryToSHE} is deferred to Section~\ref{ASHE}. 
 
\bl\label{ExpKSumLem} 
Continuing with the notations of Theorem~\ref{SHEMoment} and Proposition~\ref{AiryToSHE}, one has 
\begin{align}
\mathbb{E}_{\mathrm{Airy}}&\left[h_{k}\big(\exp(C\mathbf{a}_1, C\mathbf{a}_2, \ldots )\big)\right] = \sum_{\substack{\lambda\vdash k\\ \lambda = 1^{m_1}2^{m_2}\ldots }} \frac{\exp( \frac{k^3 T}{24})}{m_1! m_2!\ldots } \int^{\mathbf{i}\infty}_{-\mathbf{i}\infty}\frac{dw_1}{2\pi \mathbf{i}}\cdots \int^{\mathbf{i}\infty}_{-\mathbf{i}\infty} \frac{dw_{\ell(\lambda)}}{2\pi \mathbf{i}} \\\times &\mathrm{det}\left[\frac{1}{w_i+\lambda_i-w_j}\right]^{\ell(\lambda)}_{i,j=1} 
 \times \prod_{j=1}^{\ell(\lambda)} \exp\left(T/2\left[w^2_j+ (w_j+1)^2+\ldots +(w_j+\lambda_j-1)^2\right]\right).\label{eq:CIntFor}
\end{align} 
\el
\begin{proof}
Result of this lemma is implicitly present in the proof of Theorem~1.2 of \cite{BorGor16}. To avoid repetition, we present here a sketch of the proof and refer \cite{BorGor16} for more details. To begin, we define the Laplace transform of the Airy kernel as follows. For any vector $c=(c_1, \ldots ,c_n)$, we define   
\begin{align}\label{eq:LaplCor1}
R(c_1, \ldots ,c_n) = \int^{\infty}_{-\infty} e^{\sum_{i=1}^n c_ix_i} \mathrm{det}[K_{\mathrm{Ai}}(x_i,x_j)]^{n}_{i,j=1} \prod_{i=1}^{n}dx_i.
\end{align}   
One may now note that 
\begin{align}\label{eq:LaplCor2}
\mathbb{E}[h_{k}(\exp(C\mathbf{a}_1),\exp(C\mathbf{a}_2), \ldots )] = \sum_{\substack{\lambda\vdash k\\ \lambda= 1^{m_1}2^{m_2}\ldots }} \frac{1}{m_1!m_2!\ldots } R(C\lambda_1, \ldots , C\lambda_{\ell(\lambda)}).  
\end{align}  
To complete the proof, it suffices to show that the contour integral formula on the right hand side of \eqref{eq:CIntFor} is same as $R(C\lambda_1, C\lambda_2, \ldots )$. Using the following identity (see \cite[Lemma~2.6]{Oku02})
\begin{align}
\int^{\infty}_{-\infty} e^{xz} \mathrm{Ai}(z+a) \mathrm{Ai}(z+b) dz = \frac{1}{2\sqrt{\pi x}} \exp\left(\frac{x^3}{12} - \frac{a+b}{2}x- \frac{(a-b)^2}{4x}\right)
\end{align}  
we write  
\begin{align}
R(c_1, \ldots ,c_n) = \frac{e^{\sum_{i=1}^n \frac{c^3_i}{12}}}{(2\pi)^n} \int_{z_1\in \RR}dz_1\ldots \int_{z_n\in \RR}dz_n &\exp\left(-\sum_{i=1}^{n}c_iz^2_i\right)\\
&\times\mathrm{det}\left[\frac{1}{(-\mathbf{i}z_i+\frac{c_i}{2})+ (\mathbf{i}z_j +\frac{c_j}{2})}\right]^{n}_{i,j=1}. \label{eq:RForm1}
\end{align}
To show that $R(C\lambda_1, \ldots ,C\lambda_{\ell(\lambda)})$ equals to the contour integral on the right hand side of \eqref{eq:CIntFor}, we set $c_i=C\lambda_i$ and substitute $\mathbf{i}z = Cw_j+C\frac{\lambda_j}{2} - \frac{C}{2}$ for $1\leq i\leq \ell(\lambda)$in \eqref{eq:RForm}. This entails to write    
\begin{align}
R(C\lambda_n, \ldots ,C\lambda_{\ell(\lambda)}) =& \frac{e^{C^3\sum_{i=1}^{\ell(\lambda)} \frac{\lambda^3_i}{12}}}{(2\pi)^n} \int^{\mathbf{i}\infty}_{-\mathbf{i}\infty}dz_1\ldots \int^{\mathbf{i}\infty}_{-\mathbf{i}\infty}dz_{\ell(\lambda)} \\
&\times\exp\left(C^3\sum_{i=1}^{\ell(\lambda)}\lambda_i\left(w_i+\frac{\lambda_i}{2}-\frac{1}{2}\right)^2\right)\mathrm{det}\left[\frac{1}{w_j+\lambda_j-w_i}\right]^{\ell(\lambda)}_{i,j=1}. \label{eq:RForm}
\end{align}
 Applying the following identity 
 \begin{align}\label{eq:Simple1}
 \exp\left(C^3\sum_{i=1}^{\ell(\lambda)}\left(\frac{\lambda^3_i}{12}+\lambda_i\Big(w_i+\frac{\lambda_i}{2}-\frac{1}{2}\Big)^2\right)\right)
 =\exp\left(C^3\sum_{i=1}^{\ell(\lambda)}\left(w_i^2+\ldots +(w_i+\lambda_i-1)^2\right)+\frac{(Ck)^3}{12}\right)
\end{align}  
into the right hand side of \eqref{eq:RForm} indeed shows that $R(C\lambda_1, \ldots ,C\lambda_{\ell(\lambda)})$ coincides with the contour integral formula on the right hand side of \eqref{eq:CIntFor}. This completes the proof.
\end{proof}
\medskip 

\textsc{Final step in the Proof of Theorem~\ref{SHEMoment}:} Proof of Theorem~\ref{SHEMoment} follows by combining the identity \eqref{eq:kMomEq} with the formula of $\mathbb{E}_{\mathrm{Airy}}\left[h_{k}\big(\exp(C\mathbf{a}_1, C\mathbf{a}_2, \ldots )\big)\right]$ in Lemma~\ref{ExpKSumLem}.

\subsection{Proof of Proposition~\ref{AiryToSHE}}\label{ASHE}

 To prove Proposition~\ref{AiryToSHE}, we need two new inputs which are listed in Proposition~\ref{InfProdDeiv} and Proposition~\ref{SHEMomentExtract}. 

 \bp\label{InfProdDeiv}
 Continuing with the notations of Proposition~\ref{AiryToSHE}, we have 
 \begin{align}\label{eq:AiryDeri}
 \lim_{u\downarrow 0}\frac{1}{u^k}\mathbb{E}_{\mathrm{Airy}}&\left[\prod_{p=1}^{\infty}\frac{1}{1+u\exp(C\mathbf{a}_{p})}-\sum_{m=0}^{k-1} h_{m}(\exp(C\mathbf{a}_1), \exp(C\mathbf{a}_2), \ldots )(-u)^{m}\right] \\& = (-1)^{k} \mathbb{E}_{\mathrm{Airy}}\Big[h_k\big(\exp(C\mathbf{a}_1), \exp(C\mathbf{a}_2), \ldots \big)\Big] 
 \end{align} 
 for any $k\in \NN$.  
 \ep
 Proof of Proposition~\ref{InfProdDeiv} contains the core of the technical part of this paper. For clarity, we defer its proof to Section~\ref{LongProp}.

\bp\label{SHEMomentExtract} 
Recall that $\mathcal{Z}(T,X)$ is the unique solution of the SHE started from the delta initial data. Then, we have 
\begin{align}\label{eq:sheMGFDeri}
\lim_{u \downarrow 0} \frac{1}{u^k}\mathbb{E}_{\mathrm{SHE}}&\left[\exp\Big(- u\mathcal{Z}(T,0)\exp(\frac{T}{24})\Big)- \sum_{\ell=0}^{k-1}(-u)^{t}\frac{[\mathcal{Z}(T,0)]^{\ell}\exp(\frac{\ell T}{24})}{\ell!}\right] \\& = (-1)^{k}\mathbb{E}_{\mathrm{SHE}}\Big[\frac{(\mathcal{Z}(T,0))^{k}\exp(\frac{kT}{24})}{k!}\Big].   
\end{align}
\ep

\begin{proof}
Applying mean value theorem, we note
\begin{align}\label{eq:MVTappl}
\left|\exp\big(-u \theta\big) - \sum_{\ell=0}^{k} (-u)^{\ell} \frac{\theta^{\ell}}{\ell!}\right|\leq u^{k}\frac{\theta^{k+1}}{(k+1)!} \exp(-\bar{u}\theta), \quad \text{ for some }\bar{u}= \bar{u}(\theta)\in (0,u).
\end{align}  
Owing to \eqref{eq:MVTappl}, we get 
\begin{align}\label{eq:ExpUnifBd}
\frac{1}{u^{k}} \mathbb{E}&\left|\exp\big(-u \mathcal{Z}(T,0)\exp(\frac{T}{24})\big) - \sum_{\ell=0}^{k} (-u)^{\ell} \frac{(\mathcal{Z}(T,0))^{\ell}\exp(\frac{\ell T}{24})}{\ell!}\right|\\&\leq u\mathbb{E}\left[\frac{(\mathcal{Z}(T,0))^{k+1}}{(k+1)!}\exp\big(\frac{(k+1)T}{24}\big)\right]
\end{align}
where the right hand side of \eqref{eq:ExpUnifBd} is finite (see \cite[Example~2.10]{ChD15} for its upper bound). Now, letting $u$ to go to $0$ on both sides of \eqref{eq:ExpUnifBd}, we arrive at \eqref{eq:sheMGFDeri}. This completes the proof. 
\end{proof}

\subsubsection{Final step in the proof of Proposition~\ref{AiryToSHE}}
 
 We prove \eqref{eq:kMomEq} by induction. Let us first prove \eqref{eq:kMomEq} for $k=1$. Subtracting $1$ from sides of identity \eqref{eq:identity}, diving both sides by $u$, letting $u$ to go to $0$ and applying Proposition~\ref{InfProdDeiv} and Proposition~\ref{SHEMomentExtract} yield \eqref{eq:kMomEq} for $k=1$. Now, we assume that \eqref{eq:kMomEq} holds for all $k\in \ZZ_{[1,k_0]}$ where $k_0\in \NN$. To prove \eqref{eq:kMomEq} for $k=k_0+1$, using \eqref{eq:identity} and \eqref{eq:kMomEq} for $k=1,2, \ldots , k_0$, we note 
 \begin{align}\label{eq:RecentIdenty}
 \frac{1}{u^{k_0+1}}\mathbb{E}_{\mathrm{Airy}}&\left[\prod_{\ell=1}^{\infty}\frac{1}{1+u\exp(C\mathbf{a}_{\ell})}-\sum_{m=0}^{k_0} h_{m}(\exp(C\mathbf{a}_1), \exp(C\mathbf{a}_2), \ldots )(-u)^{m}\right]\\ & = \frac{1}{u^{k_0+1}}\mathbb{E}_{\mathrm{SHE}}\left[\exp\Big(- u\mathcal{Z}(T,0)\exp(\frac{T}{24})\Big)- \sum_{\ell=0}^{k_0}(-u)^{\ell}\frac{[\mathcal{Z}(T,0)]^\ell\exp(\frac{\ell T}{24})}{\ell!}\right]. 
 \end{align}
 Now, taking $u$ to $0$ on both sides of \eqref{eq:RecentIdenty} and applying Proposition~\ref{InfProdDeiv} and Proposition~\ref{SHEMomentExtract}, we arrive at \eqref{eq:kMomEq} for $k=k_0+1$. This completes the proof. 

 \subsubsection{Proof of Proposition~\ref{InfProdDeiv}}\label{LongProp}

We first prove \eqref{eq:AiryDeri} for $k=1$. Fixing $n\in \NN$, we notice  
\begin{align}\label{eq:Expand}
\prod_{p=1}^{n}\frac{1}{1+u\exp(C\mathbf{a}_{p})}-1 = -\sum_{p=1}^{n} \frac{u\exp(C\mathbf{a}_{p})}{1+u\exp(C\mathbf{a}_{p})}\prod_{m=p+1}^{n}\frac{1}{1+u\exp(C\mathbf{a}_{m})}
\end{align}
which shows that 
\begin{align}\label{eq:BaseIneq}
\left|\prod_{p=1}^{n}\frac{1}{1+u\exp(C\mathbf{a}_{p})}-1\right| \leq  u\sum_{p=1}^{n} \exp(C\mathbf{a}_{p}).
\end{align}
Taking $n$ to $\infty$ on both sides of \eqref{eq:BaseIneq}, we get 
\begin{align}\label{eq:InftyIneq}
\left|\prod_{p=1}^{\infty}\frac{1}{1+u\exp(C\mathbf{a}_{p})}-1\right| \leq  u\sum_{p=1}^{\infty} \exp(C\mathbf{a}_{p})= uh_1(\exp(C\mathbf{a}_1), \exp(C\mathbf{a}_2), \ldots ).
\end{align}
Since, the right hand side of \eqref{eq:InftyIneq} is in $L^{1}$, therefore, by the dominated convergence theorem and monotonocity of both sides of \eqref{eq:Expand}, we write   
 \begin{align}\label{eq:InfProdMinus1}
\prod_{p=1}^{\infty}\frac{1}{1+u\exp(C\mathbf{a}_{p})}-1 = -\sum_{p=1}^{\infty} \frac{u\exp(C\mathbf{a}_{p})}{1+u\exp(C\mathbf{a}_{p})}\prod_{m=p+1}^{\infty}\frac{1}{1+u\exp(C\mathbf{a}_{m})}. 
 \end{align}
 Furthermore, owing to \eqref{eq:InftyIneq}, one observes 
\begin{align}\label{eq:Ineq2}
\frac{1}{u}\left|\mathbb{E}\big[\prod_{p=1}^{\infty}\frac{1}{1+u\exp(C\mathbf{a}_p)}\big]-1\right|\leq \frac{1}{u}\mathbb{E}\left|\prod_{p=1}^{\infty}\frac{1}{1+u\exp(C\mathbf{a}_p)}-1\right|
\end{align}
 where the right hand side is bounded above by $\mathbb{E}\big[h_1(\exp(C\mathbf{a}_1), \exp(C\mathbf{a}_2), \ldots)\big].$ Thanks to \eqref{eq:Ineq2}, one can do the following interchange of limit and integral  
 \begin{align}\label{eq:Interchange}
 \lim_{u\to 0}\frac{1}{u}\Big[\mathbb{E}\Big[\prod_{p=1}^{\infty}\frac{1}{1+u\exp(C\mathbf{a}_{p})}\Big]-1\Big] = \mathbb{E}\Big[\lim_{u\to 0}\frac{1}{u}\Big[\prod_{p=1}^{\infty}\frac{1}{1+u\exp(C\mathbf{a}_{p})}-1\Big]\Big]
 \end{align}
 if the limit exists. To see why the limit on the right hand side of \eqref{eq:Interchange} exists,  we divide both sides of \eqref{eq:InfProdMinus1}  by $u$. As $u$ goes to $0$, owing to \eqref{eq:InftyIneq} and the dominated convergence theorem, one has  
 \begin{align}\label{eq:RHLimit}
 \lim_{u\to 0}\frac{1}{u}\left[\prod_{p=1}^{\infty}\frac{1}{1+u \exp(C\mathbf{a}_{p})}-1\right]= -\sum_{\ell=1}^{\infty}\exp(C\mathbf{a}_{p}) \quad \text{in }L^{1}. 
 \end{align}
 Plugging the limit of \eqref{eq:RHLimit} into the right hand side of \eqref{eq:Interchange}, we arrive at \eqref{eq:AiryDeri} for $k=1$.

 Now, we prove \eqref{eq:AiryDeri} for any $k\in \NN_{>1}$. For this, we first notice the following.
 
\medskip 
 
 \textbf{Claim:} For any $k\in \NN$, we have  
 \begin{align}\label{eq:Expj}
 \prod_{p=1}^{\infty} \frac{1}{1+u\exp(C\mathbf{a}_{p})} &= \prod_{p=1}^{\infty}\Big[\sum_{m=0}^{k}(-u\exp(C\mathbf{a}_{p}))^{m}+ \frac{(-u\exp(C\mathbf{a}_p))^{k+1}}{1+u\exp(C\mathbf{a}_{p})} \Big]
 \end{align}
 for all $u< \exp(-C\mathbf{a}_1)$.
 
\medskip 
 
 \textsc{Proof of Claim:} 
  Proof of this claim directly follows by noting that one can write 
  \begin{align}
  \frac{1}{1+u\exp(C\mathbf{a}_{p})}  = \sum_{m=0}^{k} (- u \exp(C\mathbf{a}_{p}))^{m} + \frac{(-u \exp(C \mathbf{a}_m))^{k+1}}{1+u \exp(C\mathbf{a}_{p})}\quad \forall k,p \geq 1
  \end{align}
 whenever $u < \exp(- C \mathbf{a}_1)$. 

\medskip

 \textbf{Claim:} Fix any $k\in \NN$. Then, for any $u\in ( 0, \min\{e^{\frac{-Tk^3}{24}}, \exp(-C\mathbf{a}_1)\})$,
 \begin{align}\label{eq:TrunPSer}
 \prod_{p=1}^{\infty}\big(\sum_{m=0}^{k} (- u\exp(C\mathbf{a}_p))^{m}\big) = \sum_{n=0}^{\infty} h^{(k)}_n(\exp(C\mathbf{a}_1), \exp(C\mathbf{a}_2), \ldots )(-u)^{t}
 \end{align}
 where $h^{(k)}_n$ is defined as
  \begin{align}\label{eq:htkdef}
  h^{(k)}_n\big(x_1, x_2, \ldots\big) := \sum_{(i_1, i_2, \ldots ,i_n)\in \mathcal{S}^{(k)}_n} x_{i_1}\cdot x_{i_2}\cdots x_{i_n}  
  \end{align}
  such that 
  \[
  \mathcal{S}^{(k)}_n = \{(i_1,i_2, \ldots ,i_n): i_1\leq i_2\leq \ldots \leq i_n, \text{ s.t. }|\{i_1, i_2, \ldots ,i_n\}\cap \{m\}|\leq k,\quad \forall m\in \NN\}.\]
 
 The right hand side of \eqref{eq:TrunPSer} is a power series absolutely convergent in $L^{1}$ for all $u\in (0,e^{\frac{-Tk^3}{24}})$. Moreover, for all $t\in \ZZ_{[0,k]}$, one has $h^{(k)}_t(x_1, x_2, \ldots )= h_t(x_1, x_2, \ldots )$.  

\medskip 
 
 \textsc{Proof of Claim:}
  Define the following power series
  \begin{align}\label{eq:DefHk}
 H_k(u):=\sum_{n=0}^{\infty} h^{(k)}_n\big(\exp(C\mathbf{a}_1), \exp(C \mathbf{a}_2), \ldots\big)u^{n}  
  \end{align}
  We first show that $H_k(u)$ is absolutely convergent in $L^{1}$ for all $u\in (0, \exp(-\frac{k^3T}{24}))$. Appealing to \eqref{eq:htkdef}, we note
  \begin{align}\label{eq:hkExp}
  \mathbb{E}\Big[h^{(k)}_n\big(\exp(C\mathbf{a}_1), \exp(C \mathbf{a}_2), \ldots\big)\Big] = \sum_{\substack{\lambda\vdash n\\\lambda= 1^{m_1} 2^{m_2}\ldots \\\lambda_1\leq k \ldots}} \frac{1}{m_1!m_2!\ldots } R(C\lambda_1, \ldots , C\lambda_{\ell(\lambda)}). 
  \end{align}
   Using \eqref{eq:RForm1} and Cauchy's determinantal formula \cite{C95}, we arrive at 
   \begin{align}\label{eq:RFormula}
   R(C\lambda_1, \ldots , C\lambda_{\ell(\lambda)}) = \frac{e^{\sum_{i=1}^{\ell(\lambda)}\frac{C\lambda^3_i}{12}}}{(2\pi)^{\frac{n}{2}}} \prod_{i=1}^{\ell(\lambda)}\frac{1}{(2C\lambda_i)^{\frac{3}{2}}} \mathbb{E}\left[\prod_{1\leq i<j\leq \ell(\lambda)}\frac{(Z_i-Z_j)^2+ \frac{1}{4}(\lambda_i-\lambda_j)^2}{(Z_i+Z_j)^2+ \frac{1}{4}(\lambda_i-\lambda_j)^2}\right] &&
   \end{align}
   where $Z_1, Z_2, \ldots , Z_{\ell(\lambda)} $ are i.i.d $N(0,1)$ random variables. Plugging  \eqref{eq:RFormula} into the right hand side of \eqref{eq:hkExp}, noting that that each term of the sum is bounded above by $\exp(\frac{nTk^3}{24})$ and applying Siegel's bound (see \cite[pp. 316-318]{Apostol76}, \cite[pp. 88-90]{Knopp70}) to the number of partitions of $n$, we arrive at 
   \begin{align}
   \text{l.h.s of \eqref{eq:hkExp}} \leq e^{\sqrt{n}+\frac{nTk^3}{24}}.
\end{align}    
This entails to write
\begin{align}\label{eq:ExpOfSumBd}
\mathbb{E}\left[\sum_{t=0}^{\infty}h^{(k)}_n\big(\exp(C\mathbf{a}_1), \exp(C\mathbf{a}_2), \ldots \big)u^n\right]\leq \sum_{n=0}^{\infty} u^{n}e^{\sqrt{n}+\frac{nTk^3}{24}}
\end{align}  
where the right hand side is finite for all $u\in (0, \exp(-\frac{Tk^3}{24}))$. This proves that $H_k(u)$ is absolutely convergent for all $u\in (0, \exp(-\frac{Tk^3}{24}))$ in $L^{1}$.   

It remains to show that the left hand side of \eqref{eq:TrunPSer} is indeed same as $H_k(u)$ for all $u\in (0, \min\{\exp(-\frac{Tk^3}{24}), \exp(-C\mathbf{a}_1)\})$. To prove this, for any $Q\in \NN$, we see that 
\begin{align}\label{eq:nTruncPSer}
\prod_{p=1}^{Q}\Big(\sum_{m=0}^{\infty}\big(-u\exp(C\mathbf{a}_p)\big)^m\Big) = \sum_{n=0}^{\infty}h_n^{(k)}\big(\exp(C\mathbf{a}_1), \exp(C \mathbf{a}_2), \ldots , \exp(C\mathbf{a}_Q)\big)(-u)^{n}.&& 
\end{align}
Thanks to $h^{(k)}_n(x_1, x_2, \ldots,x_Q)\leq h^{(k)}_n(x_1,x_2, \ldots )$ for all $Q\in \NN$ and $x_1,x_2,\ldots \geq 0$, one observes  
\begin{align}
|\text{l.h.s of \eqref{eq:nTruncPSer}}|\leq H_k(u). 
\end{align}
To this end, using the dominated convergence theorem and absolute convergence of $H_k(u)$, we obtain \eqref{eq:TrunPSer} for all $u\in (0, \min\{\exp(-\frac{Tk^3}{24}), \exp(-C\mathbf{a}_1)\})$.

 \medskip 
 
 \bl
  Continuing with the notation of Proposition~\ref{ExpKSumLem}, for any $k\in 2\NN$ and $u\in (0, \exp(- \frac{Tk^3}{3}))$, we define 
  \begin{align}
 \mathbf{X}_k:= \Big|\prod_{p=1}^{\infty} &\frac{1}{1+u\exp(C\mathbf{a}_{p})} -  \sum_{t=0}^{k}h_n(\exp(C\mathbf{a}_1), \exp(C\mathbf{a}_2), \ldots ) (-u)^{n}\Big|.
\end{align}   
Then $\mathbf{X}_k$ satisfies the following.
\bei 
\ii For $u<\exp(-C\mathbf{a}_1)$, 
\begin{align}\label{eq:L1ExpBd}
\mathbf{X}_k \leq \sum_{n=k+1}^{\infty} h^{(k)}_n(\exp(C\mathbf{a}_1), \exp(C\mathbf{a}_2), \ldots )u^{n} + u^{k+1}H_k(u)\sum_{p=1}^{\infty}\exp(C(k+1)\mathbf{a}_p) . &&&
\end{align} 
  See \eqref{eq:DefHk} for the definition of $H_k(.)$. 
\ii The right hand side of \eqref{eq:L1ExpBd} is in $L^{1}$ for all $u\in (0, \exp(-\frac{Tk^3}{3}))$. 

\ii We have 
 \begin{align}\label{eq:XkLimit}
\lim_{u\to 0}\frac{1}{u^k}\mathbb{E}_{\mathrm{Airy}}\mathbf{X}_k=0.
\end{align} 
\ee
\el
 
 \begin{proof}
 
\bei
  
\ii  Owing to \eqref{eq:Expj}, for all $u< \exp(-C\mathbf{a}_1)$ and $k\in 2\NN$, we observe 
 \begin{align}\label{eq:InfProdRep}
  \prod_{p=1}^{\infty} \frac{1}{1+u \exp(C\mathbf{a}_{p})} = \prod_{p=1}^{\infty}\left(\sum_{m=0}^{k} (-u \exp(C\mathbf{a}_{p}))^m\right)\prod_{p=1}^{\infty} \frac{1}{1+(u\exp(C\mathbf{a}_p))^{k+1}}. 
\end{align}  
In the same way as in \eqref{eq:InfProdMinus1}, one can now write
\begin{align}\label{eq:InfProdBreak}
\prod_{p=1}^{\infty}\frac{1}{1+(u\exp(C\mathbf{a}_{p}))^{k+1}} =1  - \sum_{p=1}^{\infty} \frac{(u\exp(C\mathbf{a}_{p}))^{k+1}}{1+(u\exp(C\mathbf{a}_{p}))^{k+1}} \prod_{m=p+1}^{\infty}\frac{1}{1+(u\exp(C\mathbf{a}_{m}))^{k+1}}. &&&
\end{align}
Plugging \eqref{eq:InfProdBreak} into the right hand side of \eqref{eq:InfProdRep}, applying \eqref{eq:TrunPSer} and noting that after subtracting $1$ from the right hand side of \eqref{eq:InfProdBreak}, the rest is bounded above by $u^{k+1}\sum_{p=1}^{\infty}\exp(C(k+1)\mathbf{a}_p)$, we arrive at \eqref{eq:L1ExpBd}. This completes showing \eqref{eq:L1ExpBd}.

\medskip 

\ii 

The first term on the right hand side of \eqref{eq:L1ExpBd} is bounded above by $H_k(u)$, thus, belongs to $L^{1}$ for all $u\in (0,\exp(-\frac{Tk^3}{12}))$. To complete proving $(ii)$, one now needs to show that 
\begin{align}\label{eq:Ek}
\mathbb{E}H_k(u)\sum_{p=1}^{\infty} \exp\big(C(k+1)\mathbf{a}_{p}\big)<\infty
\end{align}
for all $u\in (0, \exp(-\frac{Tk^3}{3}))$. To show this, we  observe 
\begin{align}
H_k(u)\sum_{p=1}^{\infty}\exp(C(k+1)\mathbf{a}_p)\leq \sum_{n=0}^{\infty}h^{(2k+1)}_{n+k+1}(\exp(C\mathbf{a}_1), \exp(C\mathbf{a}_2), \ldots )u^{n}.\label{eq:CSIineqApp}
\end{align}
which holds because of the following inequalities \[h^{(k)}_{n}(x_1,x_2, \ldots )\times \sum_{p=1}^{\infty}x^{k+1}_p\leq h^{(2k+1)}_{n+k+1}(x_1,x_2, \ldots ),\quad n\in\ZZ_{\geq 0}, k\in \NN, x_1,x_2,\ldots \in\RR_{+}\]
 Expectation of the right hand side of \eqref{eq:CSIineqApp} is finite for all $u\in (0, \exp(-\frac{T(2k)^3}{24}))$. This implies \eqref{eq:Ek}.

\ii  Write
\begin{align}\label{eq:SplitX}
\frac{1}{u^{k}}\mathbf{E}\mathbf{X}_k = \frac{1}{u^k} \mathbb{E}\mathbf{X}_k \mathbbm{1}(u<\exp(-C\mathbf{a}_1)) + \frac{1}{u^k} \mathbb{E}\mathbf{X}_k \mathbbm{1}(u\geq\exp(-C\mathbf{a}_1)). 
\end{align}
Thanks to \eqref{eq:L1ExpBd} and part $(ii)$ of this proposition, the first term on the right hand side of \eqref{eq:SplitX} converges to $0$. Note that $\mathbf{X}_k\in L^2$ and furthermore, $\mathbb{E}\mathbf{X}^2_k$ is bounded by a constant for all $u\in(0, \exp(-\frac{Tk^3}{24}))$. Applying Cauchy-Schwarz inequality, one sees 
\begin{align}\label{eq:CSIneq}
\mathbb{E}\mathbf{X}_k \mathbbm{1}(u\geq\exp(-C\mathbf{a}_1))\leq \sqrt{\mathbb{E}[\mathbf{X}^2_k]} \mathbb{P}(\mathbf{a}_1\geq -C^{-1}\log u)^{\frac{1}{2}}.
\end{align} 
The right hand side of \eqref{eq:CSIneq} is bounded above by $C\exp(-\frac{4}{3}T^{-1/2}(\log u^{-1})^{3/2})$ (thanks for instance to \cite[Theorem~1.3]{RRV11}) which if divided by $u^{k}$ converges to $0$ as $u\downarrow 0$. This shows that the second term on the right hand side of \eqref{eq:SplitX} is also converging to $0$ as $u$ nears $0$. This completes proving \eqref{eq:XkLimit}. 
\ee

\end{proof}

\medskip

\textsc{Final step in the proof of Proposition~\ref{InfProdDeiv}:} From \eqref{eq:XkLimit}, it immediately follows that
\begin{align}
\lim_{u\to 0}\frac{1}{u^k} \mathbb{E}&\left[\prod_{p=1}^{\infty} \frac{1}{1+u\exp(C\mathbf{a}_{p})} -  \sum_{n=0}^{k-1}h_n(\exp(C\mathbf{a}_1), \exp(C\mathbf{a}_2), \ldots ) (-u)^{n}\right] \\&= \mathbb{E}\Big[(-1)^{k}h_k(\exp(C\mathbf{a}_1), \exp(C\mathbf{a}_2), \ldots )\Big], \quad k\in 2\NN. 
\end{align}

Therefore, to complete the proof, it suffices to show \eqref{eq:XkLimit} holds when $k$ is odd. However, using triangle inequality, one may note that 
\begin{align}\label{eq:kodd}
\mathbf{X}_k\leq \mathbf{X}_{k+1}+ u^{k+1}h_{k+1}(\exp(C\mathbf{a}_1), \exp(C \mathbf{a}_2), \ldots).  
\end{align}  
Dividing both sides of \eqref{eq:kodd} by $u^{k}$ and letting $u$ to go to $0$, we see that right hand side of \eqref{eq:kodd} converges to $0$ in $L^{1}$ (thanks to \eqref{eq:XkLimit} for $k+1$). This proves \eqref{eq:XkLimit} when $k$ is odd.

\section{Alternative Method}\label{AlrMethod}

 In this section, we give a proof of Theorem~\ref{SHEMoment} following a suggestion of Prof. Vadim Gorin \cite{VGPrivate}. This approach is based on the computation of the limiting moment formulas of the semi-discrete directed polymers. In what follows, we recall the definition of the semi-discrete directed random polymers and its moment formulas. Subsequently, we finish proving Theorem~\ref{MainTheo} in Section~\ref{Last}.

\subsection{Semi-discrete directed random polymers}

\bd
A path of length $N$ in $\RR\times \ZZ$ is a sequence of $N$ points from $\RR\times \ZZ$. We call a path up-right and increasing if it either proceeds to the right or jumps by one unit above. For each sequence $0<s_1<\ldots <s_{N-1}<t$, we associate a path $\phi$ from $(0,1)$to $(t, N)$ such that $\phi$ jumps between the points $(s_{i},i)$ and $(s_{i+1}, i)$ for $i=1, \ldots , N-1$ and remain continuous otherwise. Let $B_1, \ldots ,B_N$ are $N$ independent Brownian motions. We define energy of the path $\phi$ as 
\begin{align}
E(\phi) = B_1(s_1) + \big[B_2(s_2) - B_2(s_1)\big]+ \ldots + \big[B_{N}(t) - B_{N}(s_{N-1})\big]. 
\end{align}   
Then, the \emph{semi-discrete polymer measure} is defined as follows 
\begin{align}
\mathbf{Z}^{N}(t) = \int e^{E(\phi)} d\phi.  
\end{align}
where the integral is defined with respect to the Lebesgue measure on the set of all up-right and increasing paths.
\ed

\bp[Proposition~5.2.8 of \cite{BC14}]\label{MomProp} 
For any $k\geq 1$, one has 
\begin{align}
\mathbb{E}\Big[ (\mathbf{Z}^{N}(t))^k\Big] = \frac{1}{(2\pi \mathbf{i})^{k}} \oint \ldots \oint \prod_{1\leq A< B \leq k} \frac{w_{A} - w_{B}}{w_{A} - w_{B}-1}\prod_{j=1}^{k}z^{-N}_j e^{t w_j} dw_j
\end{align}
where the $w_A$ contour contains only the poles at $\{w_{A+1}+1, \ldots ,w_{k}+1, 0\}$. 
\ep


In \cite{Nica16}, it was shown that the semi-discrete directed polymer measure converges to the solution of the SHE started from delta initial measure. To illustrate, let us define 
\begin{align}
C(N,T,X) = \exp\Big(N+ \frac{\sqrt{NT}+ X}{2} + X T^{-\frac{1}{2}}N^{\frac{1}{2}}\Big) (T^{\frac{1}{2}} N^{-\frac{1}{2}})^{N}.
\end{align}
\bp[Corollary~1.6 of \cite{Nica16}]\label{CDRPlim}
Fix $T\in \RR_{>0}$ and $X\in \RR$. Then,  
\begin{align}\label{eq:CDRPLim}
\mathcal{Z}^{N}(T,X) := \frac{\mathbf{Z}^{N}\big(\sqrt{TN}+X\big)}{C(N,T,X)} \Rightarrow \mathcal{Z}(T,X) \quad \text{as } N \to \infty
\end{align} 
where $\mathcal{Z}(T, X)$ is the unique solution of \eqref{eq:SHEDef} and '$\Rightarrow$' denotes the convergence in distribution.
\ep

The following proposition shows that the moments of $\mathcal{Z}^{N}(T,X)$ converge as $N$ goes to $\infty$. 

\bp
Fix $T>0$, $X\in \RR$. Then, for $k\geq 1$, we have 
\begin{align}\label{eq:MomLim}
\lim_{N\to \infty}  \mathbb{E}[(\mathcal{Z}^{N}(T,X))^k] = \frac{1}{(2\pi \mathbf{i})^{k}} \int \ldots \int \prod_{1\leq A<B\leq k} \frac{z_A-z_{B}}{z_{A} -z_{B}-1} \prod_{j=1}^{k} e^{\frac{T}{2}z^2_j + Xz_j} dz_j. 
\end{align} 
where the $z_{A}$-contour is along $C_{A}+\mathbf{i}\RR$ where $C_1>C_2+1>C_3+2>\ldots >C_k+(k-1)$. 
\ep

\begin{proof}
 Set $t = \sqrt{TN}+X$. Owing to Proposition~\ref{MomProp}, we see 
\begin{align}\label{eq:PreMomForm}
\mathbb{E}[(\mathbf{Z}^{N}(t))^k]= \frac{1}{(2\pi \mathbf{i})^{k}} \oint \ldots \oint \prod_{1\leq A< B \leq k}\frac{z_{A} - z_{B}}{z_{A} - z_{B}-1} \prod^{k}_{j=1}z^{-N}_{j} e^{t z_{j}} dz_j
\end{align} 
where the contours contain all the singularities of the integrand. Let us examine the factor $e^{tz}/ z^{N}$. Note that $tz - N\log z$ has a critical point at $z_c= N/t$. When $N$ is large, then,   
\begin{align}\label{eq:CritForm}
z_c = N/(\sqrt{TN}+X) = N^{\frac{1}{2}} T^{-\frac{1}{2}} (1- X/\sqrt{TN}).
\end{align}
Letting $z= z_c+w-\frac{X}{T}$ and Taylor expanding $e^{tz -N\log z}$ around the point $z_c$, we get 
\begin{align}
e^{tz- N\log z} &= e^{N +\sqrt{TN}(w-X/T)+ N \log (\sqrt{N/T}) -w\sqrt{TN} + O(1/N^{1/2})} e^{Xw+\frac{1}{2}T w^2}\\
& = C(N,T, X)e^{Xw+\frac{1}{2}Tw^2 + O(1/\sqrt{N})}. \label{eq:MainIntgLimit}
\end{align}
 One gets $\mathbb{E}[(\mathcal{Z}^{N}(T,X))^k]$ by dividing the right hand side of \eqref{eq:PreMomForm} by $C(N,T,X)^k$. This cancels out $C(N,T,X)^k$ which are obtained by factoring $e^{tz_j- N\log z_j}$ (see \eqref{eq:MainIntgLimit}) for $j=1, \ldots, k$. Substituting $z_j= z_c+w_j-X/T$, we see
 \begin{align}\label{eq:IntActLimit}
 \prod_{1\leq A< B\leq k}\frac{z_{A} - z_{B}}{z_{A} -z_{B}-1} = \prod_{1\leq A<B\leq k} \frac{w_{A} -w_{B}}{w_{A} - w_{B}-1}. 
 \end{align}
 Due to the exponential decay of the integrand in the region far away from the critical point, the contours transform to a sequence of bi-infinite lines which are ordered among themselves. Combining this with \eqref{eq:MainIntgLimit} and \eqref{eq:IntActLimit} and applying those to \eqref{eq:PreMomForm}, we get \eqref{eq:MomLim}.    
\end{proof}

\medskip

\subsection{Final steps of the proof of Theorem~\ref{MainTheo}}\label{Last}
Note that the limiting value in \eqref{eq:MomLim} matches with the right hand side of \eqref{eq:SHEMom}. However, it is not enough to show that the moments of the SHE are finite. As we have pointed out in Section~\ref{Intro}, one also need some uniform bound on the tail probability. To achieve this, using Markov's inequality, we see  
\begin{align}\label{eq:Markov}
\mathbb{P}\Big(\mathcal{Z}^{N}(T,X)\geq M\Big) &\leq \frac{1}{M^{k+1}} \mathbb{E}\big[(\mathcal{Z}^{N}(T,X))^{k+1}\big].  
\end{align}  
Thanks to \eqref{eq:MomLim}, we see that the right hand side of \eqref{eq:Markov} is bounded above by $C/M^{k+1}$ for all large $M$ where $C$ does not depend on $N$. Taking $N$ to $\infty$ on both sides of \eqref{eq:Markov} and combining with \eqref{eq:CDRPLim}, one obtains 
\begin{align}
\mathbb{P}\Big(\mathcal{Z}(T,X)\geq M\Big)\leq \frac{C}{M^{k+1}}, \quad \forall M>0
\end{align}  
which shows $\mathbb{E}[(\mathcal{Z}(T,X))^k]<\infty$. Owing to the uniform integrability of $\{\mathcal{Z}^{N}(T,X)\}_{N\in \NN}$  and the dominated convergence theorem,  the left hand side of \eqref{eq:MomLim} is equal to $\mathbb{E}[(\mathcal{Z}(T,X))^k]$. This completes the proof.

\section*{Acknowledgment}
This work would not have been possible without numerous inputs and valuable suggestions of Prof. Ivan Corwin  to whom the author likes to express his earnest gratitude. The author is thankful to Guillaume Barraquand for several insightful discussions and suggestions about the proof techniques. The author also likes to thank Prof. Alexei Borodin and Prof. Vadim Gorin for their comments on the early version of this paper and pointing out other methods to derive the main result of this paper.

\bibliographystyle{alpha}
\bibliography{Reference}

\end{document}